\documentclass[12pt]{amsart}
\usepackage{amsfonts}
\usepackage{amsmath,amssymb,amsthm}

\newcommand{\R}{\mathbb R}

\newcommand{\ds}{\displaystyle}
\newtheorem{thm}{Theorem}[section]

\newtheorem{lem}[thm]{Lemma}
\newtheorem{prop}[thm]{Proposition}
\theoremstyle{definition}

\theoremstyle{remark}

\newcommand\topstrut{\rule{0mm}{2.9ex}}
\newcommand\bottomstrut{\rule[-1.5ex]{0mm}{1.5ex}}
\newcommand\titlestrut{\topstrut\bottomstrut}

\begin{document}

\title[NATURAL PDE'S OF WEINGARTEN SURFACES IN EUCLIDEAN SPACE]
{Natural PDE's of linear fractional Weingarten surfaces in Euclidean
Space}
\author{Georgi Ganchev and Vesselka Mihova}%
\address{Bulgarian Academy of Sciences, Institute of Mathematics and Informatics,
Acad. G. Bonchev Str. bl. 8, 1113 Sofia, Bulgaria}%
\email{ganchev@math.bas.bg}%
\address{Faculty of Mathematics and Informatics, University of Sofia,
J. Bouchier Str. 5, 1164 Sofia, Bulgaria}
\email{mihova@fmi.uni-sofia.bg}


\begin{abstract}

We prove that the natural principal parameters on a given Weingarten
surface are also natural principal parameters for the parallel surfaces
of the given one. As a consequence of this result we obtain that the
natural PDE of any Weingarten surface is the natural PDE of its parallel
surfaces. We show that the linear fractional
Weingarten surfaces are exactly the surfaces satisfying a linear relation
between their three curvatures. Our main result is classification of
the natural PDE's of Weingarten surfaces with linear relation
between their  curvatures.

\end{abstract}

\maketitle


MS Classification: 53A05, 53A10

\section{Introduction}

The relationship between the solutions of certain types of partial differential equations
and the determination of various kinds of surfaces of constant curvature has generated
many results which have applications to the areas of both pure and applied mathematics.
This includes the determination of surfaces of either constant mean curvature or
Gaussian curvature. It has long been known that there is a connection between surfaces
of negative constant Gaussian curvature in Euclidean $R^3$ and the sine-Gordon equation.
The fundamental equations of surface theory are found to yield a type of geometrically
based Lax pair. For instance, given a particular solution of the sinh-Laplace equation,
this Lax pair can be integrated to determine the three fundamental vector fields
related to the surface.
These are also used to determine the coordinate vector field of the surface.

In \cite{Ei} Eisenhart considers three particular systems of lines on surfaces.
He founds that the so called distance-function and the radii of normal and geodesic
curvature of the directions of the lines of all three systems are simple functions
of the radii of principal curvature of the surfaces.

Further results are obtained based on the fundamental equations of surface theory,
and it is shown how specific solutions of this sinh-Laplace equation can be used
to obtain the coordinates of a surface in either Minkowski $R^3_1$
or Euclidean $R^3$ space \cite{Hu1, Hu2}.

In \cite{Br} Bracken introduces some fundamental concepts and equations pertaining
to the theory
of surfaces in three-space, and, in particular, studies a class of sinh-Laplace equation
which has the form
$\Delta u=\pm \sinh u.$

A surface $S$ with principal curvatures $\nu_1$ and $\nu_2$ is a Weingarten
surface (W-surface) \cite{W1, W2} if there exists a function $\nu$ on $S$ and two functions
(Weingarten functions) $f, \, g$ of one variable, such that
$$\nu_1=f(\nu), \quad \nu_2=g(\nu).$$

We proved in \cite {GM} that any W-surface admits locally special principal parameters
- \emph{natural principal parameters}. With respect to these natural principal parameters
the functions
$\displaystyle{ \sqrt E \exp\left(\int \frac{f'd\nu}{f-g}\right),\;
\sqrt G \exp \left(\int \frac{g'd\nu}{g-f}\right)}$ are constants,
$E,\, G$ being the coefficients of the first fundamental form on a W-surface.

With respect to natural principal parameters
any W-surface $S$ with Weingarten functions $f, \, g$ is determined uniquely up to
motions by the geometric function $\nu$, which satisfies a non-linear partial differential
equation - the \emph{natural PDE} of the surface $S$ \cite{GM}. This result solves the
Lund-Regge reduction problem \cite {FG, Fo, LR, Sym} for W-surfaces in Euclidean space.

In Proposition 3.1 we prove that

{\it The natural principal parameters of a given W-surface $S$ are natural
principal parameters for all surfaces
$\bar S(a),\; a={\rm const} \neq 0$, which are parallel to $S$.}

Theorem 3.2 states that

{\it The natural PDE of a given W-surface $S$ is the natural PDE of any
surface $\bar S(a),\; a={\rm const}\neq 0$, which is parallel to $S$.}
\vskip 2mm

To motivate our investigations let us consider surfaces in Euclidean space, whose
Gauss curvature $K$ and mean curvature $H$ satisfy the linear relation
$$\delta K=\alpha H + \gamma, \qquad \alpha^2+4\gamma\delta \neq 0.\leqno(1.1)$$

In \cite {F} it was proved that all surfaces satisfying the linear relation (1.1)
are integrable.

A similar relation has been studied in the three-dimensional
Minkowski space \cite {Mil2}.

There arises the following
 question: what are the natural PDE's describing the surfaces, whose
curvatures satisfy the relation (1.1)?

Any surface $S$, whose invariants $K$ and $H$ satisfy the linear
relation (1.1) is (locally) parallel to one of the following three types of
surfaces: a minimal surface; a CMC-surface (or a surface with positive constant
Gauss curvature); a surface with negative constant Gauss curvature.

\begin{itemize}
\item The surfaces, which are parallel to minimal surfaces, are described by
the natural PDE
$$ \lambda_{xx} + \lambda_{yy} = - e^{\lambda}.$$

\vskip 2mm
\item The surfaces, which are parallel to CMC-surfaces ($H = {\rm const}$), are
described by the one-parameter family of natural PDE's
$$\lambda_{xx} + \lambda_{yy} = -2|H| \sinh \lambda.$$

\noindent
Up to similarity, the surfaces, which are parallel to CMC-surfaces, are described by
the natural PDE of the surfaces with $H=1/2$\,.

\vskip 2mm
\item The surfaces, which are parallel to pseudo-spherical surfaces ($K={\rm const} < 0$),
are described by the one-parameter family of natural PDE's
$$ \lambda_{xx} - \lambda_{yy} = K^2 \sin  \lambda.$$

\noindent
Up to similarity, the surfaces, which are parallel to pseudo-spherical surfaces,
are described by the natural PDE of the surfaces with $K = -1$.
\end{itemize}

We call the surfaces with $H=0$, $H=1/2$ and $K=-1$ \emph{the basic classes of surfaces}
in the class of surfaces, determined by the relation (1.1).

Then we have:

{\it Up to similarity, the surfaces, whose curvatures satisfy the linear relation $(1.1)$,
are described by the natural PDE's of the basic surfaces.}
\vskip 2mm

In \cite{Mil1, Mil3} Milnor studies surface theory in Euclidean and Minkowski space,
considering harmonic maps and various relations between the curvatures $K, \,H$ and
$H'=\ds{\frac{\nu_1-\nu_2}{2}}$.

The geometric quantity $\;\rho_1-\rho_2,$ where
 $\rho_1:=(\nu_1)^{-1},\;\rho_2:=(\nu_2)^{-1}$  are the principal radii of
curvature on a given surface $S$, has a definite physical
meaning, being associated with the interval of Sturm \cite{St}, also
known as the astigmatic interval, or the amplitude of astigmatism.

A. Ribaucour \cite{Rib} has proved that
a necessary condition for the curvature lines
of the first and second focal surfaces (the first and second evolute surfaces)
of a given surface $S$ to  correspond to each
other resp. to conjugate parametric lines on $S$ is $\rho_1-\rho_2={\rm const}$ resp.
$\rho_1\,\rho_2={\rm const}$.
These are the only W-surfaces whose focal sheets have corresponding lines
of curvature \cite{Bi, Ei}.

Von Lilienthal  has proved in \cite{vL, vL1, vL2} that the first focal surface
 $\tilde S$ of a surface $S$ with $\rho_1-\rho_2={\rm const}$
is of constant negative Gauss curvature, and vice versa.

The involute surfaces  $\bar S(a), \, a \in\mathbb{R}$ of $\tilde S$
are parallel surfaces of $S$ with the property $\rho_1-\rho_2={\rm const}$.
This implies that the family $\bar S(a)$ are integrable surfaces
as a consequence of the integrability of $\tilde S$.

The curvatures of the above surfaces $S$ satisfy the relation $K=p\,H',\; p={\rm const}$.

On the Weingarten surfaces whose principal radii of curvature are bound by a relation
of the form $\rho_1=c\,\rho_2,\; c={\rm const},$
the characteristic lines \cite{Ei} cut at constant angle and only in this case.
Moreover, the conjugates of the mean orthogonal lines also cut under constant angle, so
that the configuration of all three systems is the same at all points of such a
surface \cite{Ei}.
The curvatures of these surfaces satisfy the relation $H=p\,H',\; p={\rm const}$.

Obviously the surfaces with $K=p\,H',$ or $ H=p\,H',\; p={\rm const}$ are not
included in the class characterized by (1.1).

These surfaces belong to the classes of  W-surfaces, defined by the following more general
linear relation
$$\delta K = \alpha H + \beta H' + \gamma, \quad \alpha, \beta, \gamma, \delta -
{\rm constants}; \quad \alpha^2-\beta^2 + 4 \gamma\delta \neq 0\leqno(1.2)$$
 between the Gauss curvature $K$, the mean curvature $H$ and the curvature $H'$.
We denote this class by $\mathfrak{K}$.

We show that the class $\mathfrak{K}$ is  the class of linear
fractional W-surfaces with respect to the principal curvatures.
Furthermore, if $S$ is a surface in $\mathfrak{K}$,
then its parallel surfaces $\bar S(a),\, a={\rm const},$ belong to $\mathfrak{K}$ too.

We determine ten basic relations with respect to the constants in (1.2) and each
of them generates a {\it basic subclass of surfaces} of $\mathfrak{K}$.
Any surface $S$, whose invariants $K$, $H$ and $H'$ satisfy the linear
relation (1.2) is (locally) parallel to one of these basic surfaces.

According to Theorem 3.2, we find the natural PDE's of all  surfaces of the
class $\mathfrak{K}$.

\vskip 2mm
It is essential to note that the natural PDE's of the linear fractional W-surfaces
are expressed by the following four operators:
$$\Delta \lambda:=\lambda_{xx}+\lambda_{yy},  \qquad
\bar \Delta \lambda:=\lambda_{xx}-\lambda_{yy};$$
$$\Delta^* \lambda:=\lambda_{xx}+(\lambda^{-1})_{yy},  \qquad
\bar \Delta^* \lambda:=\lambda_{xx}-(\lambda^{-1})_{yy}.$$
\vskip 3mm

The central theorem in this paper is the following
\vskip 2mm

{\bf Theorem A.} {\it Up to similarity, the surfaces in Euclidean space
free of umbilical points, whose curvatures $K$, $H$ and $H'$ satisfy the linear relation
$$\delta K = \alpha H + \beta H' + \gamma, \quad \alpha, \beta, \gamma, \delta -
{\rm constants}; \quad \alpha^2-\beta^2 + 4 \gamma\delta \neq 0,$$
are described by the natural PDE's of the following basic classes of surfaces}:

\begin{tabular}{|c|c|c|c|}
\hline
\titlestrut
Nr & Basic classes  &
{The geometric   }& Natural\\

& \titlestrut of surfaces
& function $\nu$  &
PDE\\
\hline
\titlestrut

1 & $H=0$ & $\nu=-e^{\lambda}$ & $\Delta \lambda = - e^{\lambda}$\\
[4.0ex]

2 & $H=\frac{1}{2}$ & $\nu=\frac{1}{2}(1-e^{\lambda})$ & $\Delta \lambda
= - \,\sinh{\lambda}$\\
[4.0ex]

3 & $H'=1$&  &
$\Delta^*(e^{\nu})=-2\,\nu\,(\nu+2)$\\
[4.0ex]

4 & $\begin{array}{c} H=\beta\,H',\\
\beta^2>1\end{array}$&& $\Delta^*(\nu^{\beta})=
-2\,\frac{\beta(\beta+1)}{(\beta-1)^2}\,\nu$\\
[4.0ex]

5 & $\begin{array}{c} H=\beta\,H',\\
\beta^2<1,\;\beta\neq 0\end{array}$&& $\bar\Delta^*(\nu^{\beta})=
-2\,\frac{\beta(\beta+1)}{(\beta-1)^2}\,\nu$\\
[4.0ex]

6 & $\begin{array}{c} H=\beta\,H'+1,\\
\beta^2>1\end{array}$&$\nu=\frac{(\beta-1)\,\lambda+2}{2}$&
\qquad $\Delta^*(\lambda^{\beta})=
-\frac{\beta\,((\beta-1)\lambda+2)((\beta+1)\lambda+2)}
{2\,(\beta-1)\,\lambda}$ \qquad \qquad \\
[4.5ex]

7 & $\begin{array}{c} H=\beta\,H'+1,\\
\beta^2<1,\; \beta\neq 0\end{array}$&$\nu=\frac{(\beta-1)\,\lambda+2}{2}$&
\quad $\bar\Delta^*(\lambda^{\beta})=
-\frac{\beta\,((\beta-1)\lambda+2)((\beta+1)\lambda+2)}
{2\,(\beta-1)\,\lambda}$ \quad \\
[4.5ex]

8& $K=-1$ & $\nu=\tan\lambda$  &$\bar\Delta \lambda = \sin \lambda$\\
[4.5ex]

9 & $ K=2\,H'$& $\nu=\frac{\lambda-4}{\lambda-2}
$&
$\Delta^*(e^{\lambda})=-2$\\
[4.5ex]

10 & $\begin{array}{c} K=\beta\,H'+\gamma,\\
\beta\neq 0,\,\gamma<0\end{array}$& $
\begin{array}{c}\nu=\lambda+\frac{\beta}{2},\\
[1mm]
\mathcal{I}=\frac{1}{\sqrt{-\gamma}}\,\arctan\frac{\lambda}{\sqrt{-\gamma}}
\end{array}$&
$\Delta^*(e^{\beta\,\mathcal{I}}) =\frac{\beta\,\gamma}{2}\,
\frac{\lambda\,\left(\beta\,\lambda+2\,\gamma\right)}{\lambda^2-\gamma}$\\
[4.0ex]
\hline
\end{tabular}
\vskip 6mm

The PDE's with numbers 1, 2 and 8 are the classical case of integrable equations,
which means that the corresponding classes of surfaces are also integrable.

The equation with number 9 is also integrable (cf \cite{BM, vL, vL1, vL2, vL3}).

The PDE's
with numbers 1, 2 and 8 exhaust those of them, which are expressed by
the operators $\Delta$ and $\bar \Delta$, while the ninth equation is
expressed by the operator $\Delta^*$.

For the remaining 6 types PDE's it is not known if they are integrable,
but all they are candidates to be investigated.

As an application we show that the W-surfaces with $\gamma_1=0$ are exactly the
rotational Weingarten surfaces. As examples we construct the rotational surfaces in the
classes (4) and (5) from the above theorem.

In \cite{GM1} we study analogous problems for space-like surfaces in Minkowski 3-space.
\section{Preliminaries}

In this section we introduce the standard denotations and formulas in the
theory of Weingarten surfaces in Euclidean space, which we use further.

Let $\R^3$ be the three dimensional Euclidean space with the standard flat
metric $\langle \, , \, \rangle$. We assume that the following orthonormal
coordinate system $Oe_1e_2e_3: \; e_i^2=1, \; e_i\,e_j=0, \, i\neq j$ is fixed
and gives the orientation of the space.

Let $S: \, z=z(u,v), \; (u,v) \in {\mathcal D}$\, be a surface in
$\R^3$ and $\nabla$ be the flat Levi-Civita connection of the metric
$\langle \, , \, \rangle$. The unit normal vector field to $S$ is denoted by $l$
and $E, F, G; \; L, M, N$ stand for the coefficients of
the first and the second fundamental forms, respectively.

We suppose that the surface has no umbilical points and the principal lines on
$S$ form a parametric net, i.e.
$$F(u,v)=M(u,v)=0, \quad (u,v) \in \mathcal D.$$
The principal curvatures $\nu_1, \nu_2$ and the principal geodesic
curvatures (geodesic curvatures of the principal lines) $\gamma_1,
\gamma_2$ of $S$ are given by
$$\nu_1=\frac{L}{E}, \quad \nu_2=\frac{N}{G}; \qquad
\gamma_1=-\frac{E_v}{2E\sqrt G}, \quad \gamma_2= \frac{G_u}{2G\sqrt E}. $$

We consider the tangential frame field $\{X, Y\}$ defined by
$$X:=\frac {z_u}{\sqrt E}, \qquad Y:=\frac{z_v}{\sqrt G}$$
and suppose that the moving frame $XYl$ is always positive oriented.

In what follows we consider surfaces with $\nu_1-\nu_2>0.$

The mean curvature and the Gauss curvature of $S$ are denoted as
usual by $H$ and $K$, respectively.
For our purposes we denote as the third curvature on $S$ the invariant
function
$$H':=\frac{\nu_1-\nu_2}{2}=\sqrt{H^2-K}.$$

The moving frame field $XYl$ satisfies the following Frenet type formulas:

$$\begin{tabular}{ll}
$\left|\begin{array}{llccc}
\nabla_{X} \,X & = &  &\gamma_1 \,Y + \nu_1 \, l,  &\\
[2mm]
\nabla_{X} Y & = -\gamma_1 \, X, & & \\
[2mm]
\nabla_{X} \, l & = - \nu_1 \, X; & & &
\end{array}\right.$ &
\quad
$\left|\begin{array}{llccc}
\nabla_{Y} \,X & = & & \gamma_2 \, Y, &\\
[2mm]
\nabla_{Y} Y & = -\gamma_2 \, X & &  & + \nu_2 \, l,\\
[2mm]
\nabla_{Y}\, l & = &  & - \nu_2 \, Y. &
\end{array}\right.$
\end{tabular}$$
\vskip 2mm

The integrability condition
$\nabla_X \nabla_Y l- \nabla _Y \nabla_X l-\nabla_{[X,Y]} l = 0$
for this system is equivalent to the Codazzi equations
$$\gamma_1=\frac{Y(\nu_1)}{\nu_1-\nu_2}=
\frac{(\nu_1)_v}{\sqrt G\,(\nu_1-\nu_2)}\,, \qquad
\gamma_2=\frac{X(\nu_2)}{\nu_1-\nu_2}= \frac{(\nu_2)_u}{\sqrt
E\,(\nu_1-\nu_2)}\,, \leqno(2.1)$$
and the integrability condition
$\nabla_X \nabla_Y Y- \nabla _Y \nabla_X Y-\nabla_{[X,Y]} Y= 0\,$ implies the Gauss equation
$$Y(\gamma_1)-X(\gamma_2)- (\gamma_1^2+\gamma_2^2) = \nu_1\nu_2= K. \leqno(2.2)$$

The Codazzi equations (2.1) imply the following equivalence
$$\gamma_1 \gamma_2 \neq 0 \; \iff \; (\nu_1)_v (\nu_2)_u \neq 0.$$

We consider two types of surfaces parameterized by principal parameters (cf \cite{GM}):

\begin{itemize}
\item  \emph{strongly regular} surfaces, determined by the condition
$$\gamma_1(u,v)\gamma_2(u,v) \neq 0, \quad (u,v) \in \mathcal D;$$

\item  surfaces, satisfying the conditions
$$\gamma_1(u,v)= 0,\quad \gamma_2(u,v)\neq 0, \quad (u,v) \in \mathcal D.$$

These surfaces are rotational surfaces \cite{GM} and their meridians are the first system of
principal lines.
\end{itemize}

Because of (2.1), the coefficients of the first fundamental form for strongly regular
surfaces can be expressed as functions of $\nu_1, \,\nu_2$,  $\gamma_1,\,\gamma_2$
as follows:
$$\sqrt E=\frac{(\nu_2)_u}{\gamma_2\,(\nu_1-\nu_2)} >0, \quad
\sqrt G=\frac{(\nu_1)_v}{\gamma_1\,(\nu_1-\nu_2)}>0. $$
\vskip 2mm

A surface $S:\; z=z(u,v),\; (u,v)\in \mathcal{D}$ is Weingarten
if there exist two differentiable functions
$f(\nu), \; g(\nu), \; f(\nu)-g(\nu) \neq 0, \; f'(\nu)g'(\nu)\neq 0, \;
\nu \in \mathcal{I}\subseteq{\R}$
such that the principal curvatures of $S$ at every point are given by
$\nu_1=f(\nu), \; \nu_2=g(\nu), \; \nu=\nu(u,v), \; (\nu_u(u,v),\nu_v(u,v)) \neq (0,0), \;
(u,v) \in \mathcal D$.

Let $S:\; z=z(u,v),\; (u,v)\in \mathcal{D}$ be a  Weingarten
surface parameterized by principal parameters. In \cite {GM} we proved that on $S$ the function
$$\lambda = \sqrt E \exp\left(\int \frac{f'd\nu}{f-g}\right)$$
does not depend on $v$, while the function
$$\mu = \sqrt G \exp \left(\int \frac{g'd\nu}{g-f}\right)$$
does not depend on $u$.

The principal parameters $(u,v)$ are  \emph{natural principal parameters} \cite{GM} if
$$\lambda(u)={\rm const}, \quad \mu(v) = \rm{const}.$$

Let $a= {\rm const}\neq 0,\;b={\rm const}\neq 0$, $(u_0,v_0)$ be a fixed point in
$\mathcal{D}$ and $\nu_0:=\nu(u_0,v_0)$.
The change of the parameters $(u,v)\in \mathcal D$ with
$(\bar u, \bar v)\in \bar{\mathcal{D}}$ by the formulas
$$\begin{array}{l}
\displaystyle{\bar u=a\int_{u_0}^u \sqrt E \exp\left(\int
\frac{f'd\nu}{f-g}\right)}
\, du\, +\overline{u}_0 , \;\bar u_0 = {\rm const},\\
[4mm]
\displaystyle{\bar v=b\int_{v_0}^v\sqrt G \exp \left(\int
\frac{g'd\nu}{g-f}\right)\,
dv}\,+\overline{v}_0,\; \bar v_0={\rm const} \end{array}$$
endows the surface $S$ with natural principal parameters $(\bar u, \bar v)$.

With respect to the natural principal parameters $(\bar u, \bar v)$ we get
$$ E=\frac{1}{\mathfrak{a}^2}\,\exp \left(-2\int_{\nu_0}^{\nu}
\frac{f'd\nu}{f-g}\right),\quad
 G=\frac{1}{\mathfrak{b}^2}\,\exp \left(-2\int_{\nu_0}^{\nu} \frac{g'd\nu}{g-f}\right)
\leqno(2.3)$$
with
$$\mathfrak{a}^2\, E(u_0,v_0)=1, \quad \mathfrak{b}^2\,  G(u_0,v_0)=1.$$

Let a  Weingarten surface $S:\; z=z(u,v),\; (u,v)\in \mathcal{D}$ be
parameterized by principal parameters $(u,v)$. These parameters are natural principal
if and only if \cite{GM}
$$\sqrt{EG}(\nu_1-\nu_2)={\rm const}\neq 0.\leqno(2.4)$$

The main theorem for Weingarten surfaces in [8, Theorems 5.8 and 5.13] is

\begin{thm}\label {T:7.3}
Given two differentiable functions $f(\nu), \, g(\nu); \; \nu \in \mathcal{I},$ \;
$f(\nu)-g(\nu)\neq 0$, \; $f'(\nu)g'(\nu)\neq 0$ and a differentiable function
$\nu(u,v), \; (u,v) \in {\mathcal D}$ satisfying the condition
$(\nu_u,\,\nu_v)\neq (0,0), \quad \nu(u,v)\in \mathcal{I}.$

Let $(u_0, v_0) \in \mathcal D, \; \nu_0=\nu(u_0, v_0)$ and
$\mathfrak{a}\neq 0, \, \mathfrak{b}\neq 0$
be two constants. If
$$\begin{array}{l}
\ds{\;\;\;\mathfrak{b}^2\exp\left(2\int_{\nu_0}^{\nu}
\frac{g'd\nu}{g-f}\right)\left[f'\nu_{vv}+
\left(f''-\frac{2f'^2}{f-g}\right)\nu^2_v\right]}\\
[4mm]
\ds{-\mathfrak{a}^2\exp\left(2\int_{\nu_0}^{\nu}
\frac{f'd\nu}{f-g}\right) \left[g'\nu_{uu} +
\left(g''-\frac{2g'^2}{g-f}\right)\nu^2_u\right]-fg(f-g)=0},
\end{array}\leqno(2.5)$$
then there exists a unique (up to a motion)  Weingarten surface\\
$S:\; z=z(u,v),$ \, $(u,v)\in \mathcal D_0 \subset \mathcal D$ with invariants
$$\begin{array}{c}
\nu_1=f(\nu), \quad \nu_2=g(\nu), \\
[2mm]
\displaystyle{\gamma_1=
\exp\left(\int_{\nu_0}^{\nu} \frac{g'd\nu}{g-f}\right)\,\frac{\mathfrak{b}f'}{f-g}\,\nu_v,
\; \gamma_2=
-\exp\left(\int_{\nu_0}^{\nu} \frac{f'd\nu}{f-g}\right)\,
\frac{\mathfrak{a}g'}{g-f}\,\nu_u.}
\end{array}$$
Furthermore, $(u,v)$ are natural principal parameters for $S$.
\end{thm}
\vskip 2mm

Hence, with respect to natural principal parameters each Weingarten surface possesses a
{\it natural PDE} (2.5).
\vskip 4mm

We show briefly  that  the condition
$\gamma_1 = 0$, i.e. $\nu = \nu(u)$, in Theorem 2.1
characterizes the class of rotational W-surfaces and the natural PDE of any rotational
W-surface reduces to an ODE.

In [8, Section 3] we described locally in a geometric (constructive) way the class of
surfaces whose first family $\mathcal F_1$ of principal lines consists of geodesics.

Let $c_2: x=x(v), \; v\in J_2,$ be a smooth regular curve in
${\mathbb{E}}^3$ parameterized by a natural parameter $v$ with
vector invariants $t(v), \, n(v), \, b(v)$,  curvature $\kappa
(v)>0$ and torsion $\tau (v)$.

A unit normal vector field $y(v)$ along the curve $c_2$ is said to
be \emph{torse-forming} \cite{Y} if $\nabla_t y=\alpha\, t$ \,for a certain
function $\alpha(v), \; v\in J_2$ on the curve $c_2$.

We consider an orthonormal pair $\{y_1(v), y_2(v)\}$ of
torse-forming normals along the curve $c_2\,$ and denote by
$\theta=\angle(n(v),y_1(v))$. The vector pair
$$\begin{array}{l}
y_1=\cos \theta \, n +\sin \theta \, b,\\
[2mm]
y_2=-\sin \theta \, n + \cos \theta \, b,\end{array}$$
is determined uniquely up to a constant angle $\theta_0$ by the
condition $\displaystyle{\theta(v)=-\int_0^v\tau(v)\,dv+\theta_0}$.

We choose $\theta_0=0$, i.e. the pair $\{y_1(v), y_2(v)\}$
satisfies the initial conditions $y_1(0)=n(0), \; y_2(0)=b(0)$.

The orthonormal frame field $t(v) y_1(v) y_2(v)$ satisfies the
Frenet type formulas
$$\begin{array}{l}
t'=\qquad \qquad \quad \kappa \, \cos \theta\, y_1 - \kappa \,
\sin \theta \, y_2,\\
[2mm]
y_1'=- \kappa \, \cos \theta \, t,\\
[2mm]
y_2'= \; \; \; \kappa \, \sin \theta \,t.\end{array}$$

For any $v\in J_2$ we consider the regular plane curve
$$c_1: z(s_1)=x(v) + \lambda(s_1)\, y_1(v)+ \mu(s_1)\, y_2(v) \quad s_1
\in J_1,$$
rigidly connected with every Cartesian coordinate
system $x(v)y_1(v)y_2(v)$. We suppose that the parameter $s_1$ is
natural ($\dot \lambda^2 + \dot \mu^2=1$) for $c_1$ and the functions
$\lambda, \, \mu$ satisfy the initial conditions
$\lambda(0)=\mu(0)=0;\;
 \dot\lambda(0)=0, \; \dot \mu(0)=1.$
 Then the (plane) curvature $\varkappa_1=\varkappa_1(s_1)> 0$ of $\,c_1\,$ completely
determines the functions $\lambda$ and $\mu$:
$$\lambda(s_1)=\int_0^{s_1}\sin\left(\int_0^{s_1}\varkappa_1(s_1) d s_1\right) d s_1,\quad
\mu(s_1)=\int_0^{s_1}\cos\left(\int_0^{s_1}\varkappa_1(s_1) d s_1\right) d s_1.$$

Now, let us consider the surface
$$S: \; Z(s_1,v)= x(v) + \lambda(s_1)\, y_1(v)+\mu(s_1)\, y_2(v);
\quad s_1 \in J_1, \; v \in J_2.\leqno(2.6)$$

Computing
$$\begin{array}{l}
Z_{s_1}=\dot \lambda \, y_1+\dot \mu \, y_2,\\
[2mm]
Z_v=[1-\kappa(\lambda \cos \theta - \mu \sin \theta)]\,t,\\
[2mm]
 Z_{s_1}\times Z_v=-[1-\kappa(\lambda \cos \theta - \mu \sin\theta)]
 (-\dot \mu \, y_1+ \dot \lambda \, y_2),
\end{array}$$
we get that the surface $S$ is smooth at the points, where
$$\lambda \cos \theta - \mu \sin \theta \neq \frac{1}{\kappa}.$$

We orientate the surface $S$ by choosing $l=-\dot
\mu \, y_1+ \dot \lambda \, y_2$, i.e. the normal to $S$ is the
plane normal to $c_1$, and put $X=Z_{s_1},\; Y=t$.

Let $\Gamma$ be the class of surfaces $(2.6)$
under the smoothness condition.
\vskip 2mm

{\it Any surface $S$ of the class $\Gamma$ has the following properties:

$1)$ the parametric lines are principal;

$2)$ the family $\mathcal F_1$ consists of geodesics.}
\vskip 2mm

The invariants of any surface from the class $\Gamma$ are:
$$\begin{array}{ll}
\nu_1=\varkappa_1(s_1)>0,\quad & \gamma_1=\tau_1=0;\\
[4mm]
\displaystyle{\nu_2=\frac{-\kappa\,(\lambda'\sin
\theta+ \mu'\cos\theta)}
{1-\kappa(\lambda \cos \theta - \mu \sin \theta)}}\,, & \displaystyle{\gamma_2
=\frac{-\kappa ( \lambda' \cos \theta -  \mu' \sin \theta)}
{1-\kappa(\lambda \cos \theta - \mu \sin \theta)}\,;}\\
[4mm]
\displaystyle{\varkappa_2^2=\frac{\kappa^2}{(1-\kappa(\lambda \cos \theta
- \mu \sin \theta))^2}}\,,&
\displaystyle{\tau_2=\frac{\gamma_2\, Y(\nu_2)-\nu_2\,Y(\gamma_2)}
{\gamma_2^2+\nu_2^2}}\,.
\end{array} \leqno(2.7)$$
\vskip 2mm

Theorem 3.2 \cite{GM} states that:

{\it Let $S$ be a surface parameterized by principal parameters. If the
family $\mathcal F_1$ of principal lines consists of regular
geodesics, then $S$ is locally part of a surface from the class
$\Gamma$.}
\vskip 2mm

Let now $S$ be a W-surface, parameterized by natural principal parameters $(u,v)$
and let its family $\mathcal F_1$ of principal lines consists of regular geodesics,
i.e. $\gamma_1= 0$.
The geometric function $\nu(u,v)$ in Theorem 2.1 depends in this case on $u$ only.
What is more, the invariants $\nu_1,\, \nu_2,\, \gamma_2$ of $S$ are also
functions of $u$ only and
$$\frac{d s_1}{d u}= \sqrt{E}=\frac{1}{\mathfrak{a}^2}\,\exp \left(-\int_{\nu_0}^{\nu}
\frac{f'(\nu)d\nu}{f(\nu)-g(\nu)}\right).$$

Since $S$ is at the same time a surface from the class $\Gamma$, then
the principal curvature $\nu_2$ in (2.7) does not depend of $v$. This is equivalent to the
conditions
$$(\kappa(v)\sin \theta(v))_v=(\kappa(v)\cos \theta(v))_v=0.$$

Hence $\kappa= {\rm const}>0,\; \theta={\rm const}=0$, $\,\tau_2=0$ and the curve
$c_2:\,x=x(v)$ is a circle.

Choosing for the curve $c_2$ the initial conditions $\kappa(0)=1$ and $b(0)=e_3$,
with respect to the cartesian coordinate system
$O e_1 e_2 e_3$ in $\mathbb{R}^3$ we get the meridian of the rotational W-surface
$S$ to be the curve
$$x_1=1-\lambda(u),\quad x_2=0,\quad x_3=\mu(u),$$
and its rotational axis is the coordinate axis $O x_3$.

The natural ODE of $S$ is
$$\ds{\mathfrak{a}^2\exp\left(2\int_{\nu_0}^{\nu}
\frac{f'd\nu}{f-g}\right) \left[g'\nu_{uu} +
\left(g''-\frac{2g'^2}{g-f}\right)\nu^2_u\right]+fg(f-g)=0}.$$
\vskip 2mm

\section{Parallel surfaces and their natural PDE's}

Let $S: \; z=z(u,v),\; (u,v)\in \mathcal{D}$ be a surface, parameterized by
principal parameters and $l(u,v)$ be the unit normal vector field of $S$.
The parallel surfaces of $S$ are given by
$$\bar S(a): \; \bar z(u,v)= z(u,v) + a\,l(u,v), \quad a={\rm const} \neq 0,
\quad (u,v)\in \mathcal{D}.\leqno(3.1)$$

We call the family $\{\bar S(a), \; a={\rm const} \neq 0\}$ the \emph{parallel
family} of $S$.

In this section we prove that a W-surface  and its parallel family
 have the same natural PDE.

Taking into account (3.1), we find
$$\bar z_{u}=(1-a\,\nu_1)\, z_u,\quad\bar z_{v}=(1-a\,\nu_2)\, z_v.\leqno(3.2)$$

Excluding the points, where $(1-a\,\nu_1)(1-a\,\nu_2)=0$, we obtain that the corresponding
unit normal vector fields $\bar l$ to $\bar S(a)$ and $l$ to $S$ satisfy the equality
$$\bar l=\varepsilon\,l,\quad{\rm where}\quad
\varepsilon:= {\rm sign}\,(1-a\,\nu_1)(1-a\,\nu_2).$$
Then the relations between the principal curvatures $\nu_1(u,v)$,
$\nu_2(u,v)$ of $S$ and $\bar \nu_1(u,v)$,  $\bar \nu_2(u,v)$ of its
parallel surface $\bar S(a)$ are
$$\bar \nu_1=\varepsilon\,\frac{\nu_1}{1-a\,\nu_1}\,, \quad \bar \nu_2=\varepsilon\,
\frac{\nu_2}{1-a\,\nu_2};\quad \nu_1=
\frac{\varepsilon\,\bar\nu_1}{1+a\,\varepsilon\,\bar\nu_1}\,,
\quad \nu_2=\frac{\varepsilon\,\bar\nu_2}{1+a\,\varepsilon\,\bar\nu_2}.\leqno (3.3)$$

Let $K,\; H,\;H'$ be the three invariants of the
surface $S$. The equalities (3.3) imply the  relations
between the invariants $\bar K$, $\bar H$ and $\bar H'$ of $\bar
S(a)$ and the corresponding invariants of $S$:
$$K=\frac{\bar K}{1+2 a\,\varepsilon \bar H + a^2 \bar K}\,, \quad
H=\frac{\varepsilon\,\bar H+a \bar K}{1+2 a\,\varepsilon \bar H + a^2 \bar K}\,,\quad
H'=\frac{\varepsilon\,\bar H'}{1+2 a\,\varepsilon\, \bar H + a^2 \bar
K}\,.\leqno (3.4)$$
\vskip 2mm

Now let $S: \; z=z(u,v),\; (u,v)\in \mathcal{D}$ be a Weingarten surface with Weingarten
functions $f(\nu)$ and $g(\nu)$. We suppose that $(u,v)$ are natural principal parameters
for $S$.
We show that $(u,v)$ are also natural principal parameters
for any parallel surface $\bar S(a)$.
\begin{prop}
The natural principal parameters $(u,v)$ of a given W-surface $S$ are natural
principal parameters for all parallel surfaces $\bar S(a),\; a={\rm const} \neq 0$ of $S$.
\end{prop}
{\it Proof:} Let $(u,v)\in \mathcal{D}$ be natural principal parameters for $S$,
$(u_0,v_0)$ be a fixed point in $\mathcal D$ and $\nu_0=\nu(u_0,v_0)$.
The coefficients $E$ and $G$ of the first fundamental form of $S$ are given by (2.3).
The corresponding coefficients $\bar E$ and $\bar G$ of $\bar S(a)$ in view of
(3.2) are
$$\bar E=(1-a\,\nu_1)^2\,E,\quad \bar G=(1-a\,\nu_2)^2\,G.\leqno(3.5)$$
Equalities (3.3) imply that $\bar S(a)$ is again a Weingarten surface with Weingarten
functions
$$\bar\nu_1(u,v)=\bar f(\nu)=\frac{\varepsilon f(\nu)}{1-af(\nu)}\,,
\quad\bar\nu_2(u,v)=\bar g(\nu)=\frac {\varepsilon g(\nu)}{1-ag(\nu)}. \leqno(3.6)$$
Using (3.6), we compute
$$\bar f-\bar g=\frac{\varepsilon(f-g)}{(1-a\,f)(1-a\,g)}\,,$$
which shows that ${\rm sign} \,(\bar f-\bar g)={\rm sign}\,(f-g)$.

Further, we denote by $f_0:=f(\nu_0), \; g_0:=g(\nu_0)$ and taking into account
(2.4) and (3.5), we compute
$$\sqrt{\bar E\,\bar G}\,(\bar f-\bar g)=\sqrt{E\, G}\,( f- g)=
{\rm const},$$
which proves the assertion. \qed
\vskip 2mm
Using the above statement, we prove the following theorem.

\begin{thm}
The natural PDE of a given W-surface $S$ is the natural PDE of any parallel
surface $\bar S(a),\; a={\rm const}\neq 0$, of $S$.
\end{thm}
{\it Proof}. We have to express the equation (2.5) in terms of the Weingarten functions
of the parallel surface $\bar S(a)$ of $S$.
Using (3.6), we compute successively
$$\begin{array}{l}
\displaystyle{\exp{\left(2\,\int_{\nu_0}^{\nu}
\frac{\dot{\bar g}\,d \nu}{\bar g-\bar f}\right)}
\left(\dot{\bar f}\,\nu_{vv}+\left(\ddot{\bar f}-\frac{2\,\dot{\bar
f}^2}{\bar f-\bar g}\right)\, \nu_v^2\right)}\\
[6mm]
\displaystyle{=\varepsilon\,\frac{(1-a\,g_0)^2}{(1-a\,f)^2(1-a\,g)^2}\,
\exp{\left(2\,\int_{\nu_0}^{\nu} \frac{g'\,d \nu}{ g-f}\right)}
\left(f'\,\nu_{vv}+\left( f''-\frac{2\,f'^2}{ f- g}\right)\,
\nu_v^2\right)}
\end{array}$$
and
$$\begin{array}{l}
\displaystyle{\exp{\left(2\,\int_{\nu_0}^{\nu}
\frac{\dot{\bar f}\,d \nu}{\bar f-\bar g}\right)}
\left(\dot{\bar g}\,\nu_{uu}+\left(\ddot{\bar g}-\frac{2\,\dot{\bar
g}^2}{\bar g-\bar f}\right)\, \nu_u^2\right)}\\
[6mm]
\displaystyle{=\varepsilon\,\frac{(1-a\,f_0)^2}{(1-a\,f)^2(1-a\,g)^2}\,
\exp{\left(2\,\int_{\nu_0}^{\nu} \frac{f'\,d \nu}{ f-g}\right)}
\left(g'\,\nu_{uu}+\left( g''-\frac{2\,g'^2}{ g- f}\right)\,
\nu_u^2\right).}
\end{array}$$

We have also
$$\bar f\,\bar g\,(\bar f-\bar
g)=\frac{\varepsilon}{(1-a\,f)^2(1-a\,g)^2}\, f\,g\,(f-g).$$

Putting
$$\bar E_0=(1-a\,\nu_{1}(u_0,v_0))^2\,E_0=\mathfrak{a}^{-2}\,(1-a\,f_0)^2
=:\bar {\mathfrak{a}}^{-2},$$
$$\bar G_0=(1-a\,\nu_{2}(u_0,v_0))^2\,G_0=\mathfrak{b}^{-2}\,(1-a\,g_0)^2
=:\bar {\mathfrak{b}}^{-2},$$
we compute the left hand side of (2.5) to be
$$\begin{array}{l}
\bar{\mathfrak{b}}^{2}\,\ds{\exp\left(2\int_{\nu_0}^{\nu} \frac{\dot{\bar
g}\,d\nu}{\bar g-\bar f}\right)\left(\dot{\bar f}\,\nu_{vv}+
\left(\ddot{\bar f}-\frac{2\,\dot{\bar f}^2}{\bar f-\bar g}\right)\nu^2_v\right)}\\
[6mm]
-\bar {\mathfrak{a}}^{2}\,\ds{\exp\left(2\int_{\nu_0}^{\nu}
\frac{\dot{\bar f}\,d\nu}{\bar f-\bar g}\right) \left(\dot{\bar g}\,\nu_{uu} +
\left(\ddot{\bar g}-\frac{2\,\dot{\bar g}^2}{\bar g-\bar f}\right)\nu^2_u\right)}\\
[6mm]
-\bar f\,\bar g(\bar f-\bar g)
\end{array}$$

$$\begin{array}{l}
={\mathfrak{b}}^{2}\,\ds{\exp\left(2\int_{\nu_0}^{\nu}
\frac{g'd\nu}{g-f}\right)\left(f'\nu_{vv}+
\left(f''-\frac{2f'^2}{f-g}\right)\nu^2_v\right)}\\
[6mm]
-{\mathfrak{a}}^{2}\,\ds{\exp\left(2\int_{\nu_0}^{\nu}
\frac{f'd\nu}{f-g}\right) \left(g'\nu_{uu} +
\left(g''-\frac{2g'^2}{g-f}\right)\nu^2_u\right)}\\
[6mm]
-fg(f-g).
\end{array}$$
\vskip 2mm

Hence, the natural PDE of $\bar S(a)$ in terms of the Weingarten functions
$\bar f(\nu)$, $\bar g(\nu)$ coincides with the natural PDE of $S$ in terms of the
Weingarten functions $f(\nu)$ and $g(\nu)$.

\qed
\vskip 2mm

\section{Surfaces whose curvatures satisfy a linear relation and proof of Theorem A}

We now consider surfaces, whose three
invariants $K$, $H$ and $H'$ satisfy a linear relation:
$$\delta K=\alpha \,H+\beta\,H'+\gamma, \quad \alpha, \beta, \gamma, \delta - {\rm constants},
\quad \alpha^2-\beta^2+4\gamma\delta\neq 0.\leqno(4.1)$$

In \cite {GM} we introduced linear fractional Weingarten surfaces as Weingarten surfaces
whose principal curvature functions $\nu_1$ and $\nu_2$ are related as follows
$$\nu_1=\frac{A\nu_2+B}{C\nu_2+D}, \quad A, B, C, D - {\rm constants}, \quad BC-AD \neq 0.
\leqno(4.2)$$

We show that the classes of surfaces with characterizing conditions (4.1) and (4.2),
respectively, coincide.

\begin{lem}
Any surface whose invariants $K, H, H'$ satisfy a linear relation $(4.1)$
is a linear fractional Weingarten surface determined by the relation $(4.2)$,
and vice versa.
\end{lem}

 The relations between the constants $\alpha, \beta, \gamma, \delta$ and $A, B, C, D$
are given by the equalities:
$$\alpha=A-D, \quad \beta = -(A+D), \quad \gamma=B,
\quad \delta=C.\leqno(4.3)$$

We denote by $\mathfrak K$ the class of all surfaces, free of umbilical points, whose
curvatures satisfy (4.1) or equivalently (4.2).

The aim of our study is to classify all natural PDE's of the surfaces from the class
$\mathfrak K$.
\vskip 2mm

The scheme of our investigations is the following:

The parallelism between two surfaces
given by (3.1) is an equivalence relation. On the other hand, Theorem 3.2 shows that
the surfaces from an equivalence class have one and the same natural PDE. Hence, it is
sufficient to find the natural PDE's of the equivalence classes. For any equivalence class,
we use a special representative, which we call \emph{a basic class}. Thus the classification
of the natural PDE's of the surfaces in the class $\mathfrak K$ reduces to the
classification of the natural PDE's of the basic classes.
\vskip 2mm

There are two important classes of linear fractional Weingarten surfaces:
the class of linear W-surfaces and the general class of linear fractional W-surfaces.
\vskip 3mm

{\bf I.} Let $S$ be a linear W-surface, i.e. $C= 0$ in (4.2).
Then the equality (4.1) gets the form
$$\alpha\,H+\beta\,H'+\gamma=0,\quad (\alpha,\gamma)\neq (0,0),\quad
\alpha^2-\beta^2\neq 0.\leqno(4.4)$$

In this case for the invariants of the parallel surface $\bar S(a)$ of $S$, because of (3.4),
 we get  the relation
$$\varepsilon\,(\alpha+2\,a\,\gamma)\,\bar H+\varepsilon\,\beta\,\bar H'+\gamma
=-a\,(\alpha+a\,\gamma)\,\bar K.\leqno (4.5)$$

\vskip 3mm
{\bf II.} Let $S$ be an essential linear fractional W-surface, i.e.
$C \neq 0 \; (C =1)$ in (4.2).
Then the equality (4.1) gets the form
$$K=\alpha\,H+\beta\,H'+\gamma.\leqno(4.6)$$

The corresponding relation between the invariants of the parallel surface
$\bar S(a)$ of $S$ is
$$\varepsilon(\alpha+2\,a\,\gamma)\,\bar H+\varepsilon\,\beta\,\bar H'+\gamma=
(1-a\,\alpha-a^2\,\gamma)\,\bar K.\leqno(4.7)$$
\vskip 4mm

{\bf Proof of Theorem A.}
\vskip 2mm

Let $S$ be a W-surface from the class $\mathfrak{K}$.

Each time choosing appropriate values for the constants $\mathfrak{a}$, $\mathfrak{b}$
and $\nu_0$ in (2.5) we get the following subclasses of W-surfaces of the class
$\mathfrak{K}$ and their natural PDE's.

\vskip 2mm
{\bf I. The class of linear W-surfaces.}

Let $\;\eta:= {\rm sign}\,(\alpha^2-\beta^2).$
In this case  we have the following subclasses:
\vskip 2mm

\begin{itemize}
\item[1)] $\alpha=0,\;\beta\neq 0,\;\gamma\neq 0$. Assuming that $\gamma=1$, the relation
(4.4) becomes
$$\beta\,H'+1=0. $$
\noindent
Choosing $\;\mathfrak{a}^2=e^{-\beta\,\nu_0}\;$ and  $\;\mathfrak{b}^2=e^{\beta\,\nu_0}\;$,
the natural PDE for these W-surfaces becomes
$$(e^{\beta\,\nu})_{vv} +(e^{-\beta\,\nu})_{uu}=
-\frac{2}{\beta}\,\nu\,(\beta\,\nu-2).$$
\noindent
Up to similarities these W-surfaces are generated by the basic class
$H'=1\;(\beta=-1)$ with the natural PDE
$$\Delta^*(e^{\nu})=-2\,\nu\,(\nu+2),$$
which is the class (3) in the statement of the theorem.
\vskip 2mm

\item[2)]
 $\displaystyle{\alpha\neq 0,\;\gamma=0}$. Assuming that $\alpha=1$, the relation
 (4.4) becomes
$$H+\beta\,H'=0.$$
\begin{itemize}
\item[2.1)] $\beta\neq 0,\;\eta=-1\; (\beta^2-1>0).$
\noindent
Choosing $\;\displaystyle{\mathfrak{b}^2\,\frac{\beta-1}{\beta+1}\,\nu_0^{-(\beta+1)}=1,\;
\mathfrak{a}^2\,\nu_0^{\beta-1}=1}$, the natural PDE for these W-surfaces becomes
$$\left(\nu^{\beta}\right)_{vv}+
\left(\nu^{-\beta}\right)_{uu}=
-2\,\frac{\beta(\beta-1)}{(\beta+1)^2}\,\nu.$$
\noindent
Putting $\;p=-\beta$, we get the natural PDE
$$\Delta^*(\nu^{p})=-2\,\frac{p(p+1)}{(p-1)^2}\,\nu$$
for the basic class (4) in the statement of the theorem.

\vskip 2mm
\item[2.2)] $\beta \neq 0,\;\eta=1\; (\beta^2-1<0)$.
\noindent
Choosing $\;\displaystyle{\mathfrak{b}^2\,\frac{\beta-1}{\beta+1}\,\nu_0^{-(\beta+1)}=-1,\;
\mathfrak{a}^2\,\nu_0^{\beta-1}=1}$, the natural PDE for these W-surfaces becomes
$$\left(\nu^{\beta}\right)_{vv}-
\left(\nu^{-\beta}\right)_{uu}=
2\,\frac{\beta(\beta-1)}{(\beta+1)^2}\,\nu.$$
\noindent
Putting $\;p=-\beta$, we get the natural PDE
$$\bar\Delta^*(\nu^{p})=-2\,\frac{p(p+1)}{(p-1)^2}\,\nu$$
for the basic class (5) in the statement of the theorem.
\vskip 2mm

\item[2.3)] $\beta = 0$. Putting  $\;\nu =-e^{\lambda}$, we get the natural
PDE for minimal surfaces (the elliptic Liouville equation)
$$\Delta \lambda + e^{\lambda}=0,$$
which is the basic class (1) in the statement of the theorem.
\end{itemize}
\vskip 4mm

\item[3)] $\alpha\neq 0,\;\beta=0,\; \gamma\neq 0.$ Assuming that $\alpha=1$,
the relation (4.4) becomes
$$H+\gamma=0.$$
Putting $\;\displaystyle{|H|\,e^{\lambda}:=H-\nu}=H'>0$,
we get the one-parameter system of natural PDE's for CMC surfaces
with $H=-\gamma$:
$$\Delta \lambda= -2\,|H|\,\sinh \lambda.$$
\noindent
Up to similarities these W-surfaces are generated by the basic class $|H|=\frac{1}{2}$
with the natural PDE
$$\Delta \lambda=-\sinh \lambda,$$
which is the class (2) in the statement of the theorem.
\vskip 3mm

\item[4)] $ \alpha\neq 0,\;\beta\neq 0,\;\gamma\neq 0.$
Assuming that $\alpha=1$ we have
$$H+\beta\,H'+\gamma=0,\quad \beta^2-1\neq 0.$$
Let $\displaystyle{\lambda:=2\,H'=\frac{-2}{\beta+1}\,(\nu+\gamma)}>0$.

\vskip 2mm
\begin{itemize}
\item[4.1)] If $\eta=-1 \; (\beta^2-1>0)$, we choose
$$\mathfrak{b}^2=\frac{\beta+1}{\beta-1}\,\left(\frac{-2}{\beta+1}
(\nu_0+\gamma)\right)^{\beta+1},\;
\mathfrak{a}^2\,=\left(\frac{-2}{\beta+1}(\nu_0+\gamma)\right)^{-(\beta-1)}.$$
The natural PDE becomes
$$\left(\lambda^{\beta}\right)_{vv}
+\left(\lambda^{-\beta}\right)_{uu}=
\frac{-\beta}{2\,(\beta+1)}\,
\frac{((\beta+1)\lambda+2\,\gamma)((\beta-1)\lambda+2\,\gamma)}{\lambda}\,.$$
\noindent
Up to similarities these W-surfaces are generated by the basic class\\
$H=p\,H'+1, \,p^2>1$ ($p=-\beta, \,\gamma = -1$) with the natural PDE
$$\Delta^*(\lambda^{p})=-\frac{p\,((p-1)\lambda+2)((p+1)\lambda+2)}
{2\,(p-1)\,\lambda},$$
which is the class (6) in the statement of the theorem.
\vskip 2mm

\item[4.2)] If $\eta=1\; (\beta^2-1<0)$, we choose
$$\mathfrak{b}^2=-\frac{\beta+1}{\beta-1}\,\left(\frac{-2}{\beta+1}
(\nu_0+\gamma)\right)^{\beta+1},\;
\mathfrak{a}^2\,=\left(\frac{-2}{\beta+1}(\nu_0+\gamma)\right)^{-(\beta-1)}.$$
The natural PDE becomes
$$\left(\lambda^{\beta}\right)_{vv}
-\left(\lambda^{-\beta}\right)_{uu}=
\frac{\beta}{2\,(\beta+1)}\,
\frac{((\beta+1)\lambda+2\,\gamma)((\beta-1)\lambda+2\,\gamma)}{\lambda}\,.$$
\noindent
Up to similarities these W-surfaces are generated by the basic class\\
$H=p\,H'+1, \,p^2<1$ ($p=-\beta, \,\gamma = -1$) with the natural PDE
$$\bar\Delta^*(\lambda^{p})=-\frac{p\,((p-1)\lambda+2)((p+1)\lambda+2)}
{2\,(p-1)\,\lambda},$$
which is the class (7) in the statement of the theorem.
\vskip 2mm
\end{itemize}
\end{itemize}
\vskip 2mm

{\bf II. The general class of linear fractional W-surfaces.}

We consider the subclasses:
\begin{itemize}
\item[5)] $\alpha = \gamma= 0, \; \beta \neq 0$. The relation (4.6) becomes
$$K=\beta H' \quad \iff \quad \rho_1-\rho_2=-2\,\beta^{-1}\,,$$
where $\rho_1=(\nu_1)^{-1},\;\rho_2=(\nu_2)^{-1}$ are the principal
radii of curvature of $S$.

\noindent
Putting $\lambda:=\displaystyle{4\,\frac{\nu-\beta}{2\,\nu-\beta}}$ and choosing
$\nu_0=\beta$, the natural PDE of these surfaces gets the form
$$\left(e^\lambda\right)_{uu} +\left(e^{-\lambda}\right)_{vv}+\frac{\beta^4}{8}=0.$$
\noindent
Up to similarities these W-surfaces are generated by the basic class
$K=2\,H'$ with the natural PDE
$$\Delta^*(e^{\lambda})=-2,$$
which is the class (9) in the statement of the theorem.

\vskip 3mm
\item [6)] $(\alpha, \gamma)\neq (0,0), \; \alpha^2+4\gamma\geq 0$.
The relation (4.7) implies that there exists a constant $a$, such
that $\gamma\, a^2+\alpha\, a-1=0$, so that the surface $\bar S(a)$, parallel to $S$,
has curvatures satisfying the relation (4.4).
Hence the natural PDE of this surface $\bar S(a)$ is one of the PDE's
in the linear case.
\vskip 3mm

\item[7)] $\alpha^2+4\,\gamma <0$. It follows that $\gamma <0$. The relation (4.7)
implies that there does not exist a constant $a$, such
that $\gamma\, a^2+\alpha\, a-1=0$, but for $\displaystyle{a=-\frac{\alpha}{2\,\gamma}}\;$
the surface $\bar S(a)$, parallel to $S$,
has curvatures satisfying a relation
$$K=\beta H'+\gamma.\leqno(4.8)$$

\vskip 2mm

\begin{itemize}
\item[7.1)] $\beta = 0$. The relation (4.8) becomes
$K=\gamma<0,$ i.e. $S$ is of constant negative sectional curvature $\gamma$.
Putting $\displaystyle{\lambda:=2\,\arctan\frac{\nu}{\sqrt{-\gamma}}}$,
we get the natural PDE of the surface $S$
$$\bar\Delta \lambda=K^2\,\sin \lambda.$$
\noindent
Up to similarities these W-surfaces are generated by the basic class
$K=-1$ with the natural PDE
$$\bar\Delta \lambda=\sin \lambda,$$
which is the class (8) in the statement of the theorem.
\vskip 2mm

\item[7.2)] $\beta \neq 0$, $\gamma< 0$.
 The natural PDE of $S$ is
$$(\exp{(\beta\,\mathcal{I})})_{uu}+
(\exp{(-\beta\,\mathcal{I})})_{vv}
=\frac{\beta\,\gamma}{2}\,
\frac{\lambda\,\left(\beta\,\lambda+2\,\gamma\right)}{\lambda^2-\gamma}\,,$$
where
$$\mathcal{I}=\frac{1}{\sqrt{-\gamma}}
\,\arctan\frac{\lambda}{\sqrt{-\gamma}},\quad
\lambda:=\nu-\frac{\beta}{2},$$ which is the class (10) in the
statement of the theorem.\qed
\end{itemize}
\end{itemize}

\vskip 4mm

Finally we note that any of the ten basic classes of W-surfaces contains an important
subclass consisting of rotational surfaces.
\vskip 2mm

As an example we describe and construct the rotational surfaces in the
basic classes (4) and (5) from Theorem A.

The principal curvatures of any surface $S$ in the class (4) or (5) satisfy the relation
$\displaystyle{\nu_1=\frac{\beta+1}{\beta-1}\,\nu_2} \,$ with
$\beta\neq 0,\pm 1$. Choosing $\nu_2=\nu(u,v)$, the natural parameter $s_1$ and the curvature
$\varkappa_1$ of any curve of the first family of principal lines of $S$ are respectively
\cite{GM}
$$\frac{d s_1}{d u}= \sqrt{E}=\nu^{-\frac{\beta+1}{2}},\qquad
\varkappa_1^2=\nu_1^2+\gamma_1^2=\left(\frac{\beta+1}{\beta-1}\right)^2\nu^2
\left[1+\frac{(\beta-1)^2}{4}\,\nu^{-(\beta+3)}\,\nu_v^2\right].$$
\vskip 2mm

Let now $S$ be a rotational surface and the meridians of $S$ be
the curves of the first family of principal lines. From Theorem 2.1 we get
$\gamma_1\equiv 0$ and $\nu=\nu(u)$. Thus any meridian of $S$ is given
by the parametric equations
$$\varkappa_1(u)=\left|\frac{\beta+1}{\beta-1}\,\nu\right|,\qquad
s_1(u)=\int_{u_0}^{u} \nu^{-\frac{\beta+1}{2}} d u.$$
The natural ODE of $S$ is
$$(\nu^{\beta})_{uu}=-2\,\frac{\beta(\beta+1)}{(\beta-1)^2}\,\nu,$$
and the natural
equation $\varkappa_1=\varkappa_1(s_1)$ of any meridian of $S$ is a
solution of the ODE
$$\varkappa_1 ''+\frac{2}{\beta+1}\,\varkappa_1^3=0,$$
where the derivatives of $\varkappa_1$ are taken with respect to the natural parameter
$s_1$.

We recall that the {\it Mylar Balloon} ({\it Latex Balloon})
 is constructed by taking two circular disks of
Mylar, sewing them along their boundaries and then inflating with
gas. The shape of the balloon, when it is fully inflated, is a
rotational surface with principal curvatures satisfying the equality
$\nu_1=2\nu_2$, i.e. the Mylar Balloon is a rotational surface from
the class (4) (e.g. \cite{M}).

Using the construction in section 2, it is easy to be shown that the surface
$S$ with
$$\lambda(u)=\int\nu^{-\frac{1+\beta}{2}}\sin\left(\frac{\beta+1}{\beta-1}
\int\nu^{\frac{1-\beta}{2}} d u\right) d u,\quad
\mu(u)=\int\nu^{-\frac{1+\beta}{2}}\cos\left(\frac{\beta+1}{\beta-1}
\int\nu^{\frac{1-\beta}{2}} d u\right) d u,$$
is a rotational surface in the basic class (4) ($\beta^2>1$), or (5) ($\beta^2<1$).

\end{document}